# Novel Zagreb Indices-Based Inequalities with Particular Regard to Semiregular and Generalized Semiregular Graphs

September 22, 2015

**Tamás Réti,  Imre Felde**

*Óbuda University, Bécsi út 96/B, H-1034 Budapest, Hungary*
E-mail: reti.tamas@bgk.uni-obuda.hu

**Abstract**

Topological relations between three degree-based invariants of a connected graph G are investigated. These are: the first and the second Zagreb indices ($M_1(G)$, $M_2(G)$) and the topological quantity defined as $F(G) = \sum d^3(v)$ where $d(v)$ is the degree of the vertex v, and the summation embraces all vertices of the underlying graph. We present novel inequalities including $M_1(G)$, $M_2(G)$ and $F(G)$, and show that in all cases equality holds if G is a regular or a semiregular graph. Additionally, the notion of so-called *weakly semiregular* graphs is introduced, they are considered as a possible generalization of traditional bidegreed semiregular graphs. Based on the use of Zagreb indices based graph irregularity indices, for purposes of fullerene stability prediction, comparative tests have been performed on a finite set of dual graphs of C40 fullerene isomers. By using the findings obtained, the traditional concept of graph irregularity characterization has been critically reevaluated.

## 1. Basic notions

All graphs considered in this study are finite simple connected graphs. For a graph G with n vertices and m edges, $V(G)=\{v_1, v_2, \ldots v_i, \ldots v_n\}$ and $E(G)$ denote the set of vertices and edges, respectively. To avoid trivialities, we always assume that $n \geq 3$. The degree of a vertex v is the number of edges incident to v, and denoted by $d(v)$ where $d(v) \geq 1$.

Let $d_i = d(v_i)$ be the degree of the ith vertex. Denote by d the corresponding n-element vector of degrees, and by $n_j$ the number of vertices of degree j, for $j \geq 1$. Let $\Delta = \Delta(G)$ and $\delta = \delta(G)$ be the maximum and the minimum degrees, respectively, of vertices of G. Moreover, denote by $\mu_i$ the average of the degrees of the vertices adjacent to $v_i \in V(G)$. An edge of G connecting vertices $v_i$ and $v_j$ is denoted by (i,j).

Using the standard terminology [1,2], let A=A(G) and D=D(G) be the adjacency matrix of G and the diagonal matrix of its vertex degree sequence. For a graph G, we denote by $\rho(G)$ the largest eigenvalue of A(G) and call it the spectral radius of G. As usual, the Laplacian matrix of G is L(G) = D-A, and the signless Laplacian matrix of G is Q(G)=D+A, for details see [3].

A graph is R-regular if all its vertices have the same degree R. A connected graph G is said to be bidegreed with degrees $\Delta$ and $\delta$ ($\Delta > \delta$) if at least one vertex of G has degree $\Delta$ and at least one vertex has degree $\delta$, and if no vertex of G has a degree different from $\Delta$ or $\delta$. A connected bidegreed bipartite graph $G(\Delta,\delta)$ is called *semiregular* if each vertex in the same part of bipartition has the same degree. By definition, a *complete split graph* CS(n,p) consists of an independent set with (n-p) vertices and a clique with p vertices, such that each vertex of the independent set is adjacent to each vertex of the clique [4].

Let $\varepsilon$ be a positive integer. Denote by $\Omega_\varepsilon$ the set of connected non-regular graphs for which and $|d_i - d_j| = \varepsilon$ holds for any (i,j) edge in G. If $G \in \Omega_\varepsilon$, then G is said to be a *weakly semiregular graph*. It follows that traditional bipartite semiregular graphs represent a subset of weakly semiregular graphs.

Based on the considerations outlined in [5,6,7], a connected non-regular graph G is called *nearly regular*, if $|d_i - d_j| \leq 1$ holds for any (i,j) edge in G. From the definition it follows that if G is a bidegreed nearly regular graph and $|d_i - d_j| = 1$ holds for all adjacent vertices $v_i$ and $v_j$ in G, then G is semiregular. Examples of such graphs are the complete bipartite graphs $K_{p,q}$ where q=p+1.

## 2. Preliminary considerations

In what follows, we list the previously known relevant results. Some of them will be needed in the subsequent sections. Recall that the first Zagreb index $M_1(G)$ and the second Zagreb index $M_2(G)$ of a graph G are defined as follows [8-11]:

$$M_1 = M_1(G) = \sum_{i=1}^{n} d_i^2 = \sum_{(i,j) \in E} (d_i + d_j)$$

$$M_2 = M_2(G) = \sum_{(i,j) \in E(G)} d_i d_j .$$

The Zagreb indices belong to the family of the most extensively studied molecular structure descriptors. The fundamental properties of Zagreb indices and the relationships between them are summarized in surveys [8-15].

Recently, a „forgotten topological index" denoted by F(G) has gained interest in the chemical graph theory [16,17,18]. It is defined as

$$F(G) = \sum_{i=1}^{n} d_i^3 = \sum_{(i,j)} (d_i^2 + d_j^2)$$

Some basic properties and the chemical applicability of F-index are detailed in [16].

**Lemma 1** [19]: Let G be a connected graph. Then

$$M_1(G) \leq \frac{M_2(G)}{\delta} + \delta m$$

$$M_1(G) \leq \frac{M_2(G)}{\Delta} + \Delta m$$

In both cases equality holds if G is regular or semiregular.

For a connected graph G the general sum-connectivity index $X_\alpha(G)$ proposed by Zhou and Trinajstić [20] is defined as

$$X_\alpha(G) = \sum_{(i,j) \in E} (d_i + d_j)^\alpha$$

where α is a real number. We can conclude that $X_\alpha(G)$ generalizes the first Zagreb index.

**Lemma 2** [20]: Let G be a connected graph G. If $0 < \alpha < 1$ then $X_\alpha(G) \leq M_1^\alpha(G) m^{1-\alpha}$, and if $\alpha < 0$ or $\alpha > 1$, then $X_\alpha(G) \geq M_1^\alpha(G) m^{1-\alpha}$. In both cases equality holds if and only if either G is regular or G is a bipartite semiregular graph. Depending on the choice of parameter α, several different types of topological indices can be generated.

**Corollary 1** If α= -1/2, as a particular case the ordinary sum-connectivity index denoted by X(G) is obtained [20]. For X(G) the following inequality yields:

$$X(G) = X_{-1/2}(G) = \sum_{(i,j) \in E} \frac{1}{\sqrt{d_i + d_j}} \geq \frac{m^{3/2}}{\sqrt{M_1}} .$$

Equality holds if and only if G is regular or semiregular.

**Corollary 2** Let α=2. It follows that

$$X_2(G) = \sum_{(i,j) \in E}(d_i + d_j)^2 = d^T Q d = F + 2M_2 \geq M_1^2/m$$

where Q is the signless Laplacian matrix of G. It is obvious that in the above formula equality holds if and only if G is regular or semiregular.

Because $(d_i + d_j) \leq 2\Delta$ for any (i,j) edge, this implies that

$$2\Delta \geq \frac{F + 2M_2}{M_1} = \frac{d^T Q d}{M_1} \geq \frac{M_1}{m}$$

with equalities if and only if G is regular.

**Corollary 3** For each edge (i,j) of G, the inequality $(d_i + d_j) \leq m+1$ is fulfilled, and equality holds if and only if G is isomorphic to the complete graph $K_3$ or G is isomorphic to the m-edge star graph $S_{m+1}$ [21]. Then from Corollary 2 one obtains that

$$1 \geq \frac{F + 2M_2}{(m+1)M_1} = \frac{d^T Q d}{(m+1)M_1}$$

and the equality holds if and only if $G=K_3$ or $G=S_{m+1}$.

If G is a triangle-free graph, then for each edge (i,j) of G the inequality $(d_i + d_j) \leq n$ holds with equality if and only if G is a complete bipartite graph [20]. This observation implies that

$$1 \geq \frac{F + 2M_2}{nM_1} = \frac{d^T Q d}{nM_1}$$

and the equality is attained if and only if G is isomorphic to a complete bipartite graph.

**Corollary 4** [11]: Let ρ(G) be the spectral radius of a connected graph G. Then as a consequence of Corollary 2 we get

$$\rho(G) \geq \frac{2M_2}{M_1} \geq \frac{M_1}{m} - \frac{F}{M_1}$$

and the equalities hold if and only if G is regular.

**Corollary 5** If α=3, it follows that

$$X_3(G) = \sum_{(i,j) \in E}(d_i + d_j)^3 \geq M_1^3/m^2$$

with equality if and only if G is regular or semiregular.

Because $(d_i + d_j) \leq 2\Delta$ for any (i,j) edge, we obtain

$$2\Delta \geq \frac{X_3(G)}{X_2(G)} = \frac{\sum_{(i,j) \in E}(d_i + d_j)^3}{F + 2M_2} \geq \frac{M_1^3}{(F + 2M_2)m^2} = \frac{M_1^3}{m^2 d^T Q d}$$

with equality if and only if G is regular.

Let $x^T = [x_1, x_2, \ldots x_i, \ldots x_n]$ be a vector whose elements are defined on the finite set of vertices. More exactly, let f be a non-negative continuous function defined on the vertex set V $=\{v_1, v_2, \ldots v_i, \ldots v_n\}$ of G and let $x_i = f(v_i)$ be non-negative real numbers for i=1,2,…,n.

**Lemma 3** [22]: Let G be a connected graph G. Then

$$\sum_{i=1}^n x_i = \sum_{i=1}^n f(v_i) = \sum_{(i,j) \in E} \left[ \frac{f(v_i)}{d(v_i)} + \frac{f(v_j)}{d(v_j)} \right] = \sum_{(i,j) \in E} \left[ \frac{x_i}{d(v_i)} + \frac{x_j}{d(v_j)} \right]. \quad (1)$$

**Lemma 4** As a particular case, let $x^T = d^T = [d_1, d_2, \ldots d_i, \ldots d_n]$. By using identity (1) it follows that

$$\sum_{(i,j) \in E} (d_i - d_j)^2 = \sum_{(i,j) \in E} (d_i^2 + d_j^2) - 2 \sum_{(i,j) \in E} d_i d_j = \sum_{i=1}^n d_i^3 - 2 \sum_{(i,j) \in E} d_i d_j = F - 2M_2 \geq 0 \quad (2)$$

with equality if and only if G is regular.

**Remark 1** The inequality (2) can be rewritten in the following equivalent form:

$$IRL(G) = F(G) - 2M_2(G) = \sum_{(i,j) \in E} (d_i - d_j)^2 = d^T L d \geq 0 \quad (3)$$

where L is the Laplacian matrix of G. Because the equality in (3) is fulfilled if and only if G is regular, this implies that $IRL(G) = d^T L d$ can be considered as a graph irregularity index.

**Remark 2** [16]: If G is an n-vertex connected graph, then inequality $|d_i - d_j| \leq n - 2$ is fulfilled for any (i,j) edge in G, and equality holds if and only if G is isomorphic to an n-vertex star graph $S_n$. This implies that

$$F - 2M_2 = d^T L d \leq m(n-2)^2$$

with equality if and only if G is the star $S_n$.

**Lemma 5** [23]: Let G be a connected graph G. Then

$$F = \sum_{i=1}^n d_i^3 \geq \frac{2m}{n} \sum_{i=1}^n d_i^2 = \frac{2m}{n} M_1 \geq \frac{8m^3}{n^2}. \quad (4)$$

with equality if and only if G is regular.

From inequalities represented by (3) and (4) the following conclusion can be drawn:

**Remark 3** Starting with (4) it is possible to establish a novel graph irregularity index IRF(G) defined as

$$IRF(G) = F(G) - \frac{2m}{n} M_1(G) \geq 0.$$

From the definition it follows that if graphs $G_A$ and $G_B$ have identical degree sequences then $IRF(G_A) = IRF(G_B)$ holds.

**Corollary 6** If G is a connected graph, then

$$\frac{M_1}{n} \leq \frac{\sum_{i=1}^{n} d_i^3}{2m} = \frac{F}{2m} \qquad \text{and} \qquad F = \sum_{i=1}^{n} d_i^3 \geq 2M_2.$$

In both formulas equality holds if and only if G is regular. In this particular case, the well-known Zagreb indices equality yields [13,14,15]:

$$\frac{M_1}{n} = \frac{\sum_{i=1}^{n} d_i^3}{2m} = \frac{M_2}{m}.$$

**Corollary 7** Let G be a connected graph. From Corollary 3 and Lemma 4 it follows that

$$\frac{nF}{2m} \geq M_1(G) \geq \frac{F + 2M_2(G)}{m+1} = \frac{d^T Q d}{m+1}$$

with equality if and only if G=K$_3$.

**Lemma 6** [24]: Let G be a connected graph. For G the following two sharp inequalities are valid:

$$\frac{M_1}{n} \geq \frac{4m^2}{n^2} \qquad \text{and} \qquad \frac{M_2}{m} \geq \frac{4m^2}{n^2}. \tag{5}$$

In both cases the equality holds if and only if G is regular.

**Lemma 7** Let G be a connected graph G. From Lemma 5 and Lemma 6, we get

$$F + 2M_2 = \sum_{(i,j) \in E} (d_i + d_j)^2 = d^T Q d \geq \frac{2m}{n} M_1 + 2M_2 = 2m \left( \frac{M_1}{n} + \frac{M_2}{m} \right) \geq \frac{16m^3}{n^2}$$

with equality if and only if G is regular.

**Lemma 8** [25]: Let G be a connected graph G. Then

$$F \leq (\Delta + \delta) M_1(G) - 2m\Delta\delta, \tag{6}$$

with equality if G is a regular or a bidegreed graph.

**Corollary 8** Let G be a connected graph. From [26] and Lemma 8 it follows that

$$\rho(G) \geq \sqrt{\frac{M_1}{n}} \geq \sqrt{\frac{F + 2m\Delta\delta}{n(\Delta + \delta)}}$$

and equality in each inequality holds if G is regular or semiregular.

Starting with the concept outlined in [27], a non-regular connected graph is called a *well-stabilized graph,* if any two non-adjacent vertices have equal degrees. Some different types of well-stabilized graphs are illustrated in Fig.1. The complete bipartite graphs, the non-regular complete multipartite graphs, the wheel graphs W$_n$ with n≥5 vertices, the complete split graphs CS(n,p) and the generalized windmill graphs [28] form the subsets of well-stabilized graphs. It is worth noting that the diameters of all graphs depicted in Fig.1 are equal to 2.

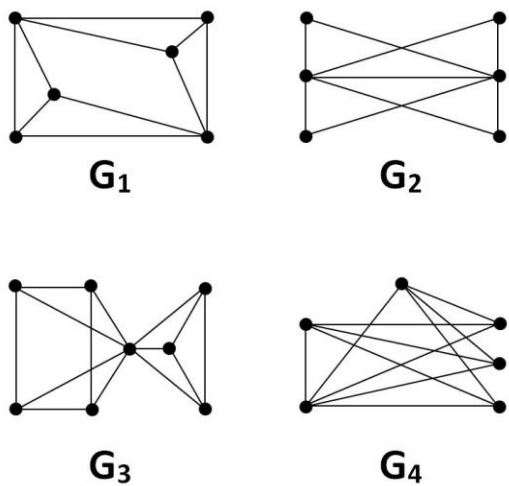

Figure 1: Examples of well-stabilized graphs

**Lemma 9** [27]: Let G be a connected graph G. Then

$$F \leq 2M_2(G) + nM_1(G) - 4m^2 \qquad (7)$$

with equality if G is a regular or a well-stabilized graph.

**Corollary 9** Consider the irregularity index VAR(G) of a connected graph G [28] defined as

$$\text{VAR}(G) = \frac{M_1(G)}{n} - \left(\frac{2m}{n}\right)^2$$

Then $F - 2M_2 = d^T L d \leq n^2 \text{VAR}(G)$ with equality if G is a regular or a well-stabilized graph.

**Corollary 10** From [11] and Lemma 9 we obtain that

$$\rho(G) \geq \frac{2M_2}{M_1} \geq \frac{F + 4m^2}{M_1} - n,$$

and equalities are attained if G is regular.

## 3. Zagreb indices based inequalities resulting in equality for semiregular graphs

In what follows, the relations between the topological indices $M_1(G)$, $M_2(G)$ and $F(G)$ will be investigated. Using these topological invariants, some novel inequalities are formulated, and it is verified that in all cases equality holds for regular and semiregular graphs.

**Proposition 1** Let G be a connected graph. Then

$$M_2(G) \geq \sum_{(i,j) \in E} \sqrt{d_i d_j \mu_i \mu_j}$$

with equality if G is regular or semiregular.

Proof. Taking into consideration the known identity formulated as

$$M_2(G) = \frac{1}{2} \sum_{i=1}^{n} d_i^2 \mu_i$$

and using the transformation rule represented by Lemma 3, we get

$$M_2(G) = \frac{1}{2} \sum_{(i,j)} (d_i \mu_i + d_j \mu_j) \geq \sum_{(i,j) \in E} \sqrt{d_i d_j \mu_i \mu_j}$$

It is obvious that equality holds for R-regular graphs, where $R = d_i = \mu_i$ holds for any vertex. If graph $G(\Delta,\delta)$ is semiregular, then $M_2(G(\Delta,\delta))=m\Delta\delta$, moreover, for an arbitrary (i,j) edge of $G(\Delta,\delta)$ the equality $d_i\mu_i = d_j\mu_j = \Delta\delta$ holds. This implies that

$$\sum_{(i,j) \in E} \sqrt{d_i d_j \mu_i \mu_j} = m\Delta\delta = M_2(G(\Delta,\delta)).$$

Several new inequalities including $M_1(G)$, $M_2(G)$ and $F(G)$ topological invariants can be easily generated by the appropriate use of the following Lemma.

**Lemma 10** Let $w_i$ be non-negative numbers (weights) for which $w_1+w_2+\ldots w_i,\ldots+ w_n =1$. Then

$$\left\{\sum_{i=1}^{n} w_i B_i\right\}^2 \leq \sum_{i=1}^{n} w_i B_i^2 \qquad (8)$$

where $B_i$ are also non-negative numbers for i=1,2,…n and the equality holds if and only if $B_1=B_2= \cdots =B_n$.

Proof. Using the Cauchy-Schwartz inequality, one gets that

$$\left\{\sum_{i=1}^{n} \sqrt{w_i}(B_i \sqrt{w_i})\right\}^2 \leq \left\{\sum_{i=1}^{n} w_i\right\}\left\{\sum_{i=1}^{n} w_i B_i^2\right\} = \sum_{i=1}^{n} w_i B_i^2 \ .$$

**Corollary 11** Let $w_i =1/n$ for i=1,2,..,n. In this particular case, the known inequality yields:

$$\left\{\frac{1}{n}\sum_{i=1}^{n} B_i\right\}^2 \leq \frac{1}{n}\sum_{i=1}^{n} B_i^2 \qquad (9)$$

**Corollary 12** Let G be an m-edge graph and denote by (i,j) an edge of G with endvertices $v_i$ and $v_j$. If B(i,j) is a non-negative number belonging to edge (i,j) in G, then

$$\left\{\frac{1}{m}\sum_{(i,j)\in E}B(i,j)\right\}^2 \leq \frac{1}{m}\sum_{(i,j)\in E}B^2(i,j) \tag{10}$$

**Corollary 13** Let G be an m-edge graph whose any edge (i,j) is characterized by a non-negative number w(i,j) for which

$$\sum_{(i,j)\in E} w(i,j) = 1.$$

If B(i,j) are also non-negative numbers, and all of them are defined as a function of edges (i,j) in G, one gets that

$$\left\{\sum_{(i,j)\in E} w(i,j)B(i,j)\right\}^2 \leq \sum_{(i,j)\in E} w(i,j)B^2(i,j). \tag{11}$$

**Proposition 2** Let G be a connected graph. Then

$$\left\{\sum_{(i,j)\in E}\sqrt{d_id_j(d_i+d_j)}\right\}^2 \leq M_1(G)M_2(G)$$

with equality if G is regular or semiregular.

Proof. Let $B(i,j) = \sqrt{d_i+d_j}/\sqrt{d_id_j}$ and $w(i,j) = d_id_j / \sum_{(i,j)} d_id_j = d_id_j / M_2(G)$.

Then using inequality (11), we have

$$\left\{\frac{1}{M_2(G)}\sum_{(i,j)\in E}\sqrt{d_id_j(d_i+d_j)}\right\}^2 \leq \frac{1}{M_2(G)}\sum_{(i,j)}(d_i+d_j) = M_1(G)/M_2(G).$$

For semiregular graphs the equality holds. This is due to the fact that if $G(\Delta,\delta)$ is a semiregular graph, then the topological quantities

$$d_i + d_j = \Delta + \delta \quad \text{and} \quad d_id_j(d_i + d_j) = \Delta\delta(\Delta+\delta)$$

are constant for any (i,j) edge in G.

Let λ be a real number, and denote by

$$^\lambda M_2(G) = \sum_{(i,j)\in E} d_i^\lambda d_j^\lambda$$

the second variable Zagreb index introduced by Miličević and Nikolić [29].

**Proposition 3** If G is a connected graph, then

$$^2M_2(G) = \sum_{(i,j)\in E} d_i^2 d_j^2 \geq \frac{M_2^2(G)}{m}$$

with equality if G is regular or semiregular.

Proof. Let $B(i,j)=d_id_j$ for any (i,j) edge in G. From inequality (10) we get

$$\frac{1}{m}\sum_{(i,j)\in E}(d_i d_j)^2 \geq \left\{\frac{1}{m}\sum_{(i,j)\in E}d_i d_j\right\}^2 = \left(\frac{M_2(G)}{m}\right)^2.$$

For a semiregular graph $G(\Delta,\delta)$ the quantity $d_i d_j = \Delta\delta$ is constant for any edge (i,j). This implies that for semiregular graphs equality holds.

**Proposition 4** If G is a connected graph, then

$$\sum_{i=1}^{n}d_i^5 + 2\sum_{(i,j)\in E}d_i^2 d_j^2 = \sum_{i=1}^{n}d_i^5 + 2(^2M_2(G)) \geq \frac{F^2}{m}.$$

Equality holds if G is regular or semiregular.

Proof. Let $B(i,j) = d_i^2 + d_j^2$. Then using inequality (10), one obtains

$$\frac{1}{m}\sum_{(i,j)\in E}(d_i^2 + d_j^2)^2 = \frac{1}{m}\left(\sum_{(i,j)\in E}(d_i^4 + d_j^4) + 2\sum_{(i,j)\in E}d_i^2 d_j^2\right) \geq \left\{\frac{1}{m}\sum_{(i,j)\in E}(d_i^2 + d_j^2)\right\}^2.$$

From this it follows that

$$m\left(\sum_{(i,j)\in E}(d_i^4 + d_j^4) + 2\sum_{(i,j)\in E}d_i^2 d_j^2\right) = m\left(\sum_{i=1}^{n}d_i^5 + 2\sum_{(i,j)\in E}d_i^2 d_j^2\right) \geq \left\{\sum_{i=1}^{n}d_i^3\right\}^2 = F^2.$$

Equality holds if $G(\Delta,\delta)$ is a semiregular graph, because for a semiregular graph $d_i^2 + d_j^2 = \Delta^2 + \delta^2$ is constant for any edge (i,j) in G.

**Proposition 5** Let G be a connected graph. Then

$$M_1(G) \geq \frac{1}{n}\left\{\sum_{i=1}^{n}\sqrt{d_i \mu_i}\right\}^2$$

with equality if G is regular or semiregular.

Proof. By definition, let $B_i = \sqrt{d_i \mu_i}$ for i= 1,2,..., n. Using inequality (9), one gets

$$\left\{\frac{1}{n}\sum_{i=1}^{n}\sqrt{d_i \mu_i}\right\}^2 \leq \frac{1}{n}\sum_{i=1}^{n}(\sqrt{d_i \mu_i})^2 = \frac{1}{n}\sum_{i=1}^{n}d_i \mu_i.$$

Using the known identity $\sum_{i=1}^{n}d_i^2 = \sum_{i=1}^{n}d_i \mu_i = M_1(G)$, it follows that

$$\left\{\frac{1}{n}\sum_{i=1}^{n}\sqrt{d_i \mu_i}\right\}^2 \leq \frac{M_1(G)}{n}.$$

The quality holds for semiregular graphs $G(\Delta,\delta)$ because $\sqrt{d_i \mu_i} = \sqrt{\Delta\delta}$ is constant for any vertex in G, moreover the identity $M_1(G(\Delta,\delta))=m(\Delta+\delta)=n\Delta\delta$ is fulfilled.

**Proposition 6** Let G be a connected graph. Then

$$\frac{M_1(G)}{n} \leq \frac{1}{M_1(G)} \sum_{i=1}^{n} (d_i \mu_i)^2 .$$

and equality holds if G is regular or semiregular.

Proof. Let $B_i = d_i \mu_i$ for i= 1,2,…,n. Using inequality (9) we have

$$\left\{ \frac{1}{n} \sum_{i=1}^{n} (d_i \mu_i) \right\}^2 = \frac{M_1^2(G)}{n^2} \leq \frac{1}{n} \sum_{i=1}^{n} (d_i \mu_i)^2 .$$

**Corollary 14** From [25] we obtain that

$$\rho^2(G) \geq \frac{\sum_{i=1}^{n} (d_i \mu_i)^2}{M_1} \geq \frac{M_1}{n},$$

and equality in each inequality is simultaneously attained if G is regular.

**Proposition 7** Let G be an n-vertex connected graph. Then

$$\sum_{i=1}^{n} \mu_j^2 \geq \frac{1}{n} \left\{ \sum_{i=1}^{n} |d_i - \mu_i| \right\}^2 + M_1(G)$$

with equality if G is regular or semiregular.

Proof. By definition, let $B_i = |d_i - \mu_i| \geq 0$ for i= 1,2,…, n. Using inequality (9), one gets

$$\left\{ \frac{1}{n} \sum_{i=1}^{n} |d_i - \mu_i| \right\}^2 \leq \frac{1}{n} \sum_{i=1}^{n} (d_i - \mu_i)^2 = \frac{1}{n} \left\{ \sum_{i=1}^{n} d_i^2 + \sum_{i=1}^{n} \mu_i^2 - 2 \sum_{i=1}^{n} d_i \mu_i \right\}.$$

Because $\sum_{i=1}^{n} d_i^2 = \sum_{i=1}^{n} d_i \mu_i = M_1(G)$, it follows that

$$\left\{ \frac{1}{n} \sum_{i=1}^{n} |d_i - \mu_i| \right\}^2 \leq \frac{1}{n} \sum_{i=1}^{n} (d_i - \mu_i)^2 = \frac{1}{n} \left\{ \sum_{i=1}^{n} \mu_i^2 - M_1(G) \right\}.$$

If $G=G(\Delta,\delta)$ is a semiregular graph, then $|d_i - \mu_i| = \Delta - \delta$ is a positive constant for any vertex $v_i$ in G. Consequently, we have

$$\left\{ \frac{1}{n} \sum_{i=1}^{n} |d_i - \mu_i| \right\}^2 = \{\Delta - \delta\}^2 = \frac{1}{n} \sum_{i=1}^{n} (d_i - \mu_i)^2 = \frac{1}{n} \sum_{i=1}^{n} (\Delta - \delta)^2 .$$

**Proposition 8** Let G be a connected graph. Then

$$\frac{1}{m} \left\{ \sum_{(i,j) \in E} |\sqrt{\mu_i} - \sqrt{\mu_j}| \right\}^2 \leq M_1(G) - 2 \sum_{(i,j) \in E} \sqrt{\mu_i \mu_j}$$

with equality if G is regular or semiregular.

Proof. Let $B(i, j) = \left|\sqrt{\mu_i} - \sqrt{\mu_j}\right|$. Starting with inequality (10) we get

$$\left\{\frac{1}{m}\sum_{(i,j)\in E}\left|\sqrt{\mu_i} - \sqrt{\mu_j}\right|\right\}^2 \leq \frac{1}{m}\left\{\sum_{(i,j)\in E}\left(\sqrt{\mu_i} - \sqrt{\mu_j}\right)^2\right\} = \frac{1}{m}\left\{\sum_{(i,j)\in E}(\mu_i + \mu_j) - 2\sum_{(i,j)\in E}\sqrt{\mu_i\mu_j}\right\}$$

Now, using the Lemma 2, we get $M_1(G) = \sum_{i=1}^n d_i\mu_i = \sum_{(i,j)\in E}(\mu_i + \mu_j)$. This implies that

$$\frac{1}{m}\left\{\sum_{(i,j)\in E}\left|\sqrt{\mu_i} - \sqrt{\mu_j}\right|\right\}^2 \leq \sum_{(i,j)\in E}\left(\sqrt{\mu_i} - \sqrt{\mu_j}\right)^2 = M_1(G) - 2\sum_{(i,j)\in E}\sqrt{\mu_i\mu_j}.$$

If $G=G(\Delta,\delta)$ is a semiregular graph, equality yields because $\left|\sqrt{\mu_i} - \sqrt{\mu_j}\right| = \sqrt{\Delta} - \sqrt{\delta}$ is a positive constant for any (i,j) edge in G.

**Proposition 9** Let G be a connected graph. Then

$$X_3(G) = \sum_{(i,j)\in E}(d_i + d_j)^3 \geq \frac{1}{M_1(G)}\{F(G) + 2M_2(G)\}^2 = \frac{\{d^TQd\}^2}{M_1(G)}.$$

In the above formula equality is fulfilled if G is regular or semiregular.

Proof. Let $B(i,j) = d_i + d_j$ and $w(i, j) = (d_i + d_j)/\sum_{(i,j)}(d_i + d_j) = (d_i + d_j)/M_1(G)$.

Starting with the inequality (11), we have

$$\frac{1}{M_1(G)}\sum_{(i,j)\in E}(d_i + d_j)(d_i + d_j)^2 \geq \left(\frac{1}{M_1(G)}\sum_{(i,j)\in E}(d_i + d_j)(d_i + d_j)\right)^2.$$

It follows that

$$\frac{1}{M_1(G)}\sum_{(i,j)\in E}(d_i + d_j)^3 \geq \left(\frac{1}{M_1(G)}\sum_{(i,j)\in E}(d_i + d_j)^2\right)^2.$$

Finally, one gets

$$M_1(G)\sum_{(i,j)\in E}(d_i + d_j)^3 \geq \left(\sum_{(i,j)\in E}(d_i + d_j)^2\right)^2 = \{F(G) + 2M_2(G)\}^2 = \{d^TQd\}^2.$$

It is easy to see that in the above formula equality is fulfilled if G is regular or semiregular.

**Corollary 15** Because $(d_i + d_j) \leq 2\Delta$ for any (i,j) edge, we obtain on the one hand

$$2\Delta M_1(G) \geq F(G) + 2M_2(G) = d^TQd,$$

and the other hand

$$\Delta - \frac{F}{2M_1} \geq \frac{M_2}{M_1}$$

with equality in both cases if G is regular.

**Proposition 10** If G is a connected graph then

$$X_{1/2}^2(G) = \left\{ \sum_{(i,j) \in E} \sqrt{d_i + d_j} \right\}^2 \leq M_2(G) \sum_{(i,j) \in E} \frac{d_i + d_j}{d_i d_j}.$$

Equality is fulfilled if G is regular or semiregular.

Proof. Let $B(i,j) = \sqrt{d_i + d_j}/d_i d_j$ and $w(i,j) = d_i d_j / \sum_{(i,j)} d_i d_j = d_i d_j / M_2(G)$. Now, using inequality (11) we have

$$\left\{ \frac{1}{M_2(G)} \sum_{(i,j) \in E} d_i d_j \frac{\sqrt{d_i + d_j}}{d_i d_j} \right\}^2 = \left\{ \frac{1}{M_2(G)} \sum_{(i,j) \in E} \sqrt{d_i + d_j} \right\}^2 \leq \frac{1}{M_2(G)} \sum_{(i,j)} d_i d_j \left\{ \frac{\sqrt{d_i + d_j}}{d_i d_j} \right\}^2.$$

For semiregular graphs the equality holds. If $G(\Delta, \delta)$ is a semiregular graph, then the topological quantities

$$d_i + d_j = \Delta + \delta \quad \text{and} \quad (d_i + d_j)/d_i d_j = (\Delta + \delta)/\Delta\delta$$

are constant for any (i,j) edge in G.

## 4. Some properties of nearly regular and weakly semiregular graphs

The bidegreed nearly regular graphs have several relevant applications in the mathematical chemistry. The benzenoid graphs [30] having vertices with degrees 2 and 3 belong to the family of nearly regular graphs. Similarly, the set of bidegreed dual graphs of classical fullerene graphs with vertex degrees 5 and 6 is considered as a special class of nearly regular graphs.

As we have already mentioned, weakly semiregular graphs are obtained as a possible generalization of bidegreed semiregular graphs. There exists a particular class of connected non-regular graphs for which $|d_i - d_j| = 1$ holds for any (i,j) edge in G. For simplicity, such graphs are called *weakly irregular graphs*. From this definition it follows that weakly irregular graphs are weakly semiregular and nearly regular simultaneously.

In Fig.2 different types of non-regular graphs are shown. Graph $H_1$ is a bidegreed weakly semiregular graph with $\varepsilon=2$, graph $H_2$ is a tridegreed weakly irregular graph, while graphs denoted by $H_3$ and $H_4$ are nearly regular graphs.

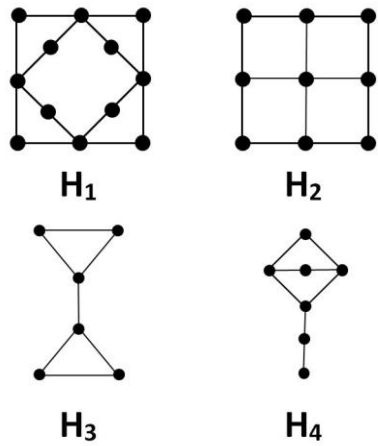

Figure 2: Some different types of non-regular graphs

Consider the Albertson irregularity index [31] of a graph G defined as

$$\text{Irr}(G) = \sum_{(i,j) \in E} |d_i - d_j|. \tag{12}$$

**Proposition 11** Let G be a connected graph. Then

$$\text{Irr}(G) \leq F(G) - 2M_2(G) = d^T L d \tag{13}$$

with equality if G is regular or nearly regular.

Proof. On the one hand, one obtains that

$$F(G) - 2M_2(G) = d^T L d = \sum_{(i,j) \in E} (d_i - d_j)^2 \geq \sum_{(i,j) \in E} |d_i - d_j| = \text{Irr}(G)$$

on the other hand, equality holds if

$$(d_i - d_j)^2 = |d_i - d_j| = 0 \quad \text{or} \quad (d_i - d_j)^2 = |d_i - d_j| = 1$$

for any edge (i,j) in G. It is easy to see that the above equalities are fulfilled for regular or nearly regular graphs.

**Corollary 16** If G is a weakly irregular graph with m edges, then

$$\text{Irr}(G) = F(G) - 2M_2(G) = d^T L d = m.$$

**Proposition 12** Let G be a connected graph. Then

$$\text{Irr}^2(G) = \left\{ \sum_{(i,j) \in E} |d_i - d_j| \right\}^2 \leq m(F(G) - 2M_2(G)) = m d^T L d.$$

In the above formula equality holds if and only if, G is regular or weakly semiregular.

Proof. Let $B(i, j) = |d_i - d_j|$. Using inequality (10) we obtain

$$\left(\frac{1}{m}\sum_{(i,j)\in E}|d_i - d_j|\right)^2 \leq \frac{1}{m}\sum_{(i,j)\in E}(d_i - d_j)^2.$$

It follows that

$$\left(\sum_{(i,j)\in E}|d_i - d_j|\right)^2 \leq m\sum_{(i,j)\in E}(d_i - d_j)^2 = m\left\{\sum_{i=1}^n d_i^3 - 2\sum_{(i,j)\in E}d_id_j\right\}.$$

It is obvious that equality holds if G is regular. It is easy to see that equality is also fulfilled for weakly semiregular graphs where $|d_i - d_j| = \varepsilon \geq 1$ holds for any (i,j) edge of G.

**Proposition 13** If G is a weakly semiregular graph with $\varepsilon \geq 1$, then G is bipartite.

Proof. It is known that (i) every acyclic graph is bipartite, and (ii) a graph is bipartite if and only if it does not contain an odd cycle [1]. Consequently it is enough to verify that if G is a weakly semiregular graph with $\varepsilon \geq 1$, then G does not contain odd cycles. Let $C_M$ be a cycle in G which includes M edges. Denote by $\{v_1, v_2,\ldots v_i, v_{i+1},\ldots , v_{M-1}, v_M\}$ the sequence of vertices in $C_M$, and by $\{d_1, d_2,\ldots d_i, d_{i+1},\ldots ,d_{M-1}, d_M\}$ the corresponding degrees of vertices, where $v_{M+1}=v_1$ and $d_{M+1}=d_1$, by definition. Because $C_M$ is a cycle in a weakly semiregular graph, this implies that $d_i \geq 2$, $d_{i+1} - d_i = +\varepsilon$ or $-\varepsilon$, and

$$\mathrm{Sum}(C_M) = (d_1 - d_2) + (d_2 - d_3) + \ldots + (d_M - d_1) = 0.$$

Denote by $M_p$ and $M_n$ the numbers of edges for which $(d_{i+1} - d_i)$ are positive and negative numbers, respectively. $\mathrm{Sum}(C_M)=0$ is fulfilled only if $M_p = M_n$. In this case $M=2M_p=2M_n$, consequently the edge number M in the cycle $C_M$ should be an even number.

**Remark 4** Let $\Delta \geq 2$ be an arbitrary integer and denote by $\{1, 2,\ldots, \Delta-1, \Delta\}$ an increasing sequence of $\Delta$ positive integers. There exists a weakly irregular tree graph $T_\Delta$ characterized unambiguously by the degree set $\{1, 2,\ldots, \Delta-1, \Delta\}$.

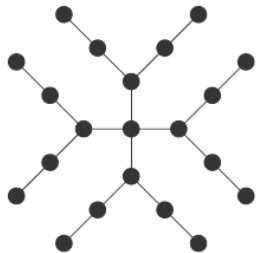

Figure 3: Dendrimer graph $T_\Delta$ (case of $\Delta=4$)

The construction of the uniquely defined tree graph $T_\Delta$ is illustrated in Fig. 3.

## 5. Chemical application: Structural characterization of fullerene isomers

By using the inequalities considered it is possible to define several molecular descriptors of new types by which the structure of molecules can be characterized and compared on the basis of quantitative criteria. In what follows, the application of some novel molecular descriptors is demonstrated on an example concerning the structural classification of all carbon fullerene molecules.

A fullerene (fullerene graph) with k vertices, denoted by $C_k$ exists for all even k≥ 20 except k=22, where the number of pentagons is 12, and the number of hexagons is k/2-10 [32]. In our investigation, for testing purpose, the set of C40 fullerene isomers was selected. As it is known, in the majority of cases, for the stability prediction of lower fullerene isomers $C_k$ with k≤70 the so-called pentagon adjacency index Np is used [32,33,34]. By definition, Np is equal to the total number of edges between adjacent pentagons. It is supposed that fullerenes which minimize $N_P$ are more likely to be stable than those that do not [32-38].

When analyzing and testing the ability of a topological invariant for prediction purposes, it is important to take into consideration the following concept: it is required that a selected graph invariant (a topological descriptor of a fullerene) must be able not only to distinguish between isomers structures but also to rank the isomers in the order of decreasing (or increasing) stability. For the structural characterization of fullerenes we used the dual graphs of the original fullerene graphs [38]. The dual graph $C_k^{dual}$ of a traditional fullerene isomer $C_k$ is a bidegreed planar graph, it contains 12 vertices of degree 5, and any remaining vertices are of degree 6. Three topological invariants defined on the set of C40 dual graphs were selected for testing purpose: the traditional pentagon index Np, and the indices IRL(G) and IRM(G). The IRL(G) index is represented by formula (3)

$$\text{IRL}(G) = F(G) - 2M_2(G) = \sum_{(i,j)\in E}(d_i - d_j)^2 = d^T L d \geq 0. \tag{14}$$

The IRM(G) index [38] is defined as

$$\text{IRM}(G) = M_2(G) - \frac{4m^3}{n^2} = m\left(\sqrt{\frac{M_2(G)}{m}} - \frac{2m}{n}\right)\left(\sqrt{\frac{M_2(G)}{m}} + \frac{2m}{n}\right) \geq 0. \tag{15}$$

Both of them belong to the family of graph irregularity indices, they are non-negative numbers, and IRM(G)=IRL(G)=0 if and only if G is a regular graph.

Starting with the previous considerations it should be emphasized that the duals graphs of fullerene graphs belong to the family of nearly regular graphs, (the inequality $|d_i - d_j| \leq 1$ holds for any (i,j) edge of $C_k^{dual}$). This implies that according to formula (13) the relation between the Albertson irregularity index Irr(G) and the IRL(G) index of $C_k^{dual}$ dual graphs reduces to the following simple identity:

$$IRL(C_k^{dual}) = d^T L d = Irr(C_k^{dual})$$

From the definition of Np pentagon index it follows [22] that for a dual fullerene graph $C_k^{dual}$:

$$IRLd = IRL(C_k^{dual}) = Irr(C_k^{dual}) = 60 - 2Np . \qquad (16)$$

Moreover, performing simple calculations it is easy to show [38] that for $C_k^{dual}$ dual graphs

$$IRMd = IRM(C_k^{dual}) = 54k + Np - 360 - \frac{54k^3}{(k+4)^2} \qquad (17)$$

As can be seen, the pentagon adjacency index is included in both formulas, and IRL(G) and IRM(G) are linear functions of Np. From this observation we can conclude that for classical fullerenes the discriminating power of IRL(G) and IRM(G) indices is considered to be equivalent.

Simultaneously, using the Density Functional Tight-Binding (DFTB) method [39] we calculated the total tight binding energy TQ characterizing the relative stability of isomers. On the basis of the calculated energy values, the relative stability of fullerene isomers can be predicted with high accuracy. It is generally supposed that the most stable isomers are characterized by the lowest values of this energy [36,37,38,39].

In Fig.4 for the forty $C_{40}$ isomers the variation of the tight binding energy TQ with IRLd and IRMd indices is shown. From diagrams depicted in Fig.4 we can conclude that IRLd and IRMd topological descriptors can be efficiently used for stability prediction of fullerene isomers, both of them have a very good correlation with the total tight binding energy. The coefficient of determination measuring the strength of linear relationship is relatively high, its values are equal ($r^2$=0.955) for both diagrams displayed in Fig. 4.

This finding can be explained by the fact, that IRLd and IRMd indices are the linear functions of pentagon adjacency index Np.

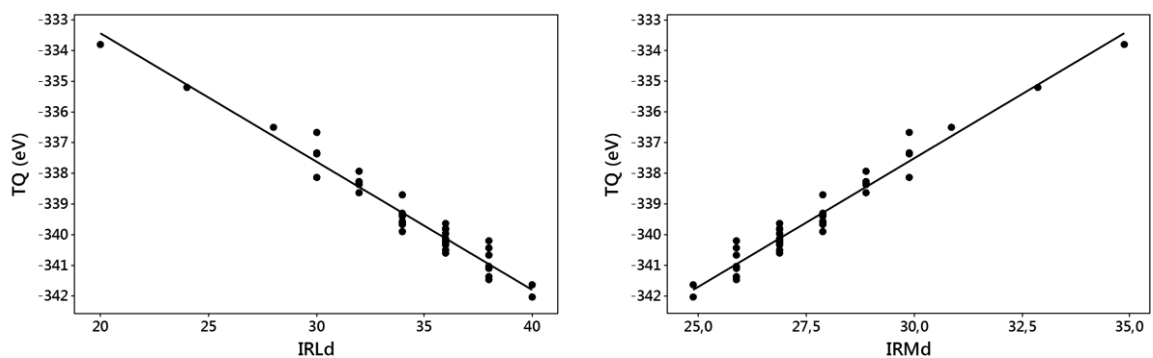

Figure 4: Relation between IRLd and IRMd indices and total tight binding energy

In Fig. 5, the correspondence between the graph irregularity indices IRLd and IRMd is demonstrated.

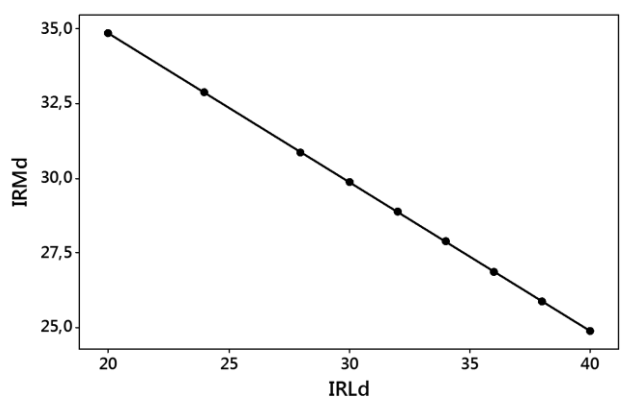

Figure 5: Relationship between topological indices IRLd and IRMd for $C_{40}$ isomers

The diagram shown in Fig.5 is strange and surprising. Comparing the computed IRLd and IRMd indices for C40 isomers, the following interesting results are obtained:

i) Eliminating the pentagon index from (16 and 17), we get a linear equation representing a simple correspondence between IRLd and IRMd indices:

$$\mathrm{IRMd} = 54k - 330 - \frac{54k^3}{(k+4)^2} - \frac{\mathrm{IRLd}}{2}$$

where k is the vertex number of an arbitrary $C_k$ isomer. From this equation it follows that IRLd and IRMd characterize the irregularity of dual graphs of fullerene isomers in a revised order.

ii) As an example, for C40:38 and C40:0 isomers the following values of computed graph irregularity indices were obtained (see Fig. 5):

24.876=IRMd(C40:38) < IRMd(C40:1)=34.876 and 40=IRLd(C40:38) > IRLd(C40:1)=20.

As can be observed, for the irregularity classification of C40:38 and C40:0 isomers we have completely different results with IRLd and IRMd indices. (It is known that among forty C40 isomers, fullerene C40:38 is the most stable, and the least stable isomer is C40:1 [32,33,34].) On the basis of our findings it can be concluded that there exist infinitely many graphs for which the irregularity ranking performed by IRLd and IRMd is not compatible with our "subjective" expectations.

It is important to emphasize that from equations (14) and (15) it follows directly that the linear relationship

$$\mathrm{IRM}(G) = \frac{1}{2}F(G) - \frac{4m^3}{n^2} - \frac{\mathrm{IRL}(G)}{2}. \qquad (18)$$

This is generally valid for all sets of n-vertex connected graphs having an identical vertex degree sequence. As an example, consider the 8-vertex chemical tree graphs $T_1$, $T_2$ and $T_3$ depicted in Fig.6 [40].

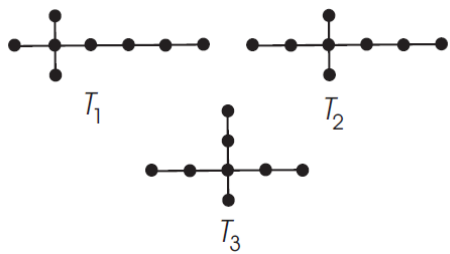

Figure 6: Structurally similar trees with identical vertex degree distribution [40]

They are characterized by an identical vertex degree sequence ($n_1$=4, $n_2$=3, $n_4$=1) and the following topological parameters: $M_2(T_1)$=30, $M_2(T_2)$=32, $M_2(T_3)$=37 and $F(T_j)$=92 for j=1,2,3. Performing simple computations, as a particular case of (18) we obtain:

$$\text{IRM}(T_A) = 24.5625 - \frac{1}{2}\text{IRL}(T_A).$$

where $T_A$ denotes any 8-vertex tree with vertex degree sequence: ($n_1$=4, $n_2$=3, $n_4$=1).

Our observations demonstrate clearly that in several cases the usual concept applied to graph irregularity characterization is false, because it is based on "ill-founded subjective considerations". This strange phenomenon called "the graph irregularity paradox", roughly speaking, can be interpreted as follows: There exist particular sets of connected graphs for which the result of evaluation (ranking) of graph irregularity by using various irregularity indices is determined primarily on their definitions, consequently the results of classification depend only to a limited extent on the true topological structure of graphs considered.

## 6. Final remarks

Formulas (10) and (11) can be efficiently used to generate novel inequalities between different degree-based topological indices. This is demonstrated on the following example. Consider the ordinary sum-connectivity index $X(G)$ and the harmonic index $H(G)$ of a graph G [20,21]:

$$X(G) = \sum_{(i,j) \in E} \frac{1}{\sqrt{d_i + d_j}} \quad , \quad H(G) = \sum_{(i,j) \in E} \frac{2}{d_i + d_j}.$$

Now, using the inequality (10) in which the non-negative numbers $B(i,j)$ are defined as $B(i,j) = 1/\sqrt{d_i + d_j}$, one obtains that

$$H(G) \geq \frac{2}{m} X^2(G) = \frac{2}{m} \left\{ \sum_{(i,j) \in E} \frac{1}{\sqrt{d_i + d_j}} \right\}^2.$$

It is easy to see that in the above formula equality holds if and only if G is regular or semiregular. According to Corollary 1, for any connected graph G the inequality $X(G) \geq m^{3/2}/M_1^{1/2}$ is fulfilled. This implies that

$$H(G) \geq \frac{2}{m} X^2(G) \geq \frac{2m^2}{M_1(G)}$$

where equality holds if and only if either G regular or G is a bipartite semiregular graph. Moreover, as a particular case, for an R-regular graph $G_R$ with n vertex and m edge, one obtains that

$$H(G_R) = \frac{2}{m} X^2(G_R) = \frac{2m^2}{M_1(G_R)} = \frac{n}{2}.$$